\documentclass[12pt]{article}

\usepackage[a4paper, top=2.54cm, bottom=2.54cm, left=2cm, right=2cm]{geometry}  
\usepackage{wallpaper}  

\usepackage{anyfontsize}

\usepackage[nottoc]{tocbibind} 
\usepackage[titletoc,title]{appendix}
\usepackage{tikz} 

\usepackage{amssymb,amsmath,amsthm}	
\usepackage{lscape} 
\usepackage{bm} 
\usepackage{enumerate}
\usepackage{tikz-cd}
\usetikzlibrary{matrix,arrows,decorations.pathmorphing}
\usepackage{mathrsfs} 
\numberwithin{equation}{section} 
\usepackage{graphicx}

\newtheorem{thm}{Theorem}[section]
\newtheorem{prop}[thm]{Proposition}

\theoremstyle{definition}
\newtheorem{defi}[thm]{Definition}
\newtheorem{ex}[thm]{Example}
\theoremstyle{remark}

\theoremstyle{plain}

\newcommand{\Z}{\mathbb {Z}}
\newcommand{\C}{\mathbb {C}}
\newcommand{\PP}{\mathbb {P}}

\newcommand{\R}{{\rm \bf R}}
\newcommand{\Q}{{\rm \bf Q}}
\newcommand{\bQ}{{\rm \bf \underline{Q}}}

\newcommand{\GL}{{\rm GL}}
\newcommand{\Id}{{\rm Id}}
\newcommand{\SL}{{\rm SL}}

\newcommand{\im}{{\rm im}}

\newcommand{\End}{{\rm End}}

\newcommand{\Cone}{{\rm Cone}}

\newcommand{\cN}{\mathcal{N}}
\newcommand{\cC}{\mathcal{C}}

\newcommand{\fg}{\mathfrak {g}}

\newcommand{\fgl}{\mathfrak {gl}}

\newcommand{\fS}{\mathfrak {S}}

\newcommand{\unit}{M_{\langle1\rangle}^{\oplus n}}

\title{Maximal Border Subrank Tensors}
\author{Chia-Yu Chang}
\date{\today}

\begin{document}
\pagenumbering{roman}
\begin{center}
    \begin{Huge}
        Maximal Border Subrank Tensors\\
    \end{Huge}
        
\vspace{30pt}

Chia-Yu Chang

\vspace{20pt}

July 30, 2022
\vspace{10pt}
\end{center}
 
\renewcommand{\abstractname}{Abstract}
\begin{abstract}
We prove a lower bound on the dimension of the set of maximal border subrank tensors. This is the first such bound of its type. 
\end{abstract}
 





\setcounter{page}{1}
\pagenumbering{arabic}

\section{Introduction}
Let $A,B,C$ be $n$ dimensional vector spaces over $\C$, and let $\{a_i\}$, $\{b_i\}$, $\{c_i\}$ respectively be bases of $A,B,C$, and $\{\alpha_i\}$, $\{\beta_i\}$, $\{\gamma_i\}$ the dual bases. Let $M_{\langle1\rangle}^{\oplus r}:=\sum_{i=1}^r a_i\otimes b_i\otimes c_i$ for any positive integer $r$. For $r=n$, we call $\unit$ the unit tensor of $A\otimes B\otimes C$.

\begin{defi}
    An element of $A\otimes B\otimes C$ is said to have {\it rank one} if it is of the form $a\otimes b\otimes c$ where $a\in A$, $b\in B$, and $c\in C$.
    Let $T$ be an element in $A\otimes B\otimes C$. The {\it rank} of $T$, denote $\R(T)$, is defined as the minimal positive integer $r$ such that $T=\sum_{i=1}^rZ_i$ with each $Z_i$ rank one.
    The {\it border rank} of $T$, denote $\underline{\R}(T)$, is defined as the minimal positive integer $r$ such that $T$ is a limit of tensors of rank $r$.
\end{defi}

\begin{defi}
    The {\it subrank} of $T$, denote $\Q(T)$, is defined as the maximal positive integer $s\leq n$ such that $$M_{\langle1\rangle}^{\oplus s}\in(\End(A)\times\End(B)\times\End(C))\cdot T.$$ The {\it border subrank} of $T$, denote $\bQ(T)$, is the maximal positive integer $s\leq n$ such that $$M_{\langle1\rangle}^{\oplus s}\in \overline{\GL(A)\times\GL(B)\times\GL(C)\cdot T}.$$ 
\end{defi}

Note that the above definitions are independent of the choice of bases. Among the four quantities, border subrank is the least understood. This paper is about border subrank.

\subsection{Motivation from complexity theory}
The complexity of the matrix multiplication tensor is measured by the {\it exponent of matrix multiplication}, a constant that is denoted $\omega$. It is defined as the smallest number such that for all $\epsilon>0$, the rank of the $n\times n\times n$ matrix multiplication tensor is $O(n^{\omega+\epsilon})$. It is known that $\omega$ is less than $3$ and no less than $2$. In 1969 \cite{MR248973}, Strassen showed the first non-trivial upper bound $\omega\leq\log_27$. Several methods for finding an upper bound of $\omega$ or trying to prove $\omega=2$ have been developed. A {\it barrier} \cite{MR3984631} for a method is any lower bound for all upper bounds on $\omega$ that can be obtained by that method. When a tensor is of big subrank, it is less likely that it will affected by the barriers, and when it is of maximal border subrank, it is not subject to any barrier. On the other hand, if the border rank is large, since it is not known how to compute the asymptotic rank, it is effectively useless. The tensor rank measures the ``cost'' of a tensor and the subrank measures the ``value'' of a tensor.

A well-know method to find upper bounds on $\omega$ is the {\it laser method}, which is based on the work of Strassen \cite{MR709378}. The idea is to study an intermediate tensor $T$ which can be proven to have low border rank and is ``close to'' being a matrix multiplication tensor. If some big Kronecker power $T^{\boxtimes N}$ restricts to a big direct sum $\oplus_iM_{\langle n_i,n_i,n_i\rangle}$ of matrix multiplication tensors, then we have an upper bound. Let $\Tilde{\R}(T):=\lim_{N\to\infty}\left(\R(T^{\boxtimes N})\right)^{1/N}$ be the asymptotic rank of $T$, and let $\Tilde{\Q}(T):=\lim_{N\to\infty}\left(\Q(T^{\boxtimes N})\right)^{1/N}$ be the asymptotic subrank of $T$. When $\log\Tilde{\R}(T)/\log\Tilde{\Q}(T)>1$, the tensor $T$ can not be used to prove that $\omega=2$.

Note that asymptotic rank and asymptotic subrank may also be defined using border rank and border subrank. Then the tensors of maximal border subrank are not subject to barriers for their utility in the laser method. More generally this is true of tensors of maximal asymptotic subrank, but that is more difficult to determine. Thus one motivation for the paper is to expand the list of tensors that could potentially be used to prove $\omega=2$ via the laser method.

The rank and border rank of a generic tensor are the same and equal to the maximal border rank. However, the behavior of the subrank and border subrank are different from the behavior of the rank and border rank. In 2022 \cite{https://doi.org/10.48550/arxiv.2205.15168}, Derksen, Makam, and Zuiddam proved that the generic subrank of tensors have bounds $3(\lfloor\sqrt{n/3+1/4}-1/2\rfloor)\leq\Q(n)\leq\lfloor\sqrt{3n-2}\rfloor.$ In particular the generic subrank is not maximal. We still do not have a good upper bound for generic border subrank and it wasn't even known what the set of maximal border subrank tensors is, neither was its dimension. We give a lower bound of the dimension of the set of maximal border subrank tensors and a description of some members in that set.

We consider the nullcone determined by the symmetry group of the unit tensor. The Hilbert-Mumford criterion, see e.g., \cite[Theorem 2.4.3]{MR1918599}, implies that the nullcone can be determined using a maximal torus of the symmetry group. The sum of the unit tensor and a general element in the nullcone is of maximal border subrank. So by counting the dimension of the set of such tensors and its orbit closure, we find a lower bound of the dimension of the set of all maximal border subrank tensors.

\subsection{Subrank and border subrank}

For any tensor $T$, it is clear that $\bQ(T)\geq\Q(T)$.
By the definition of subrank, the maximal possible subrank and border subrank for a tensor in $A\otimes B\otimes C$ is $n$ and the unit tensor $\unit$ has the maximal subrank. The following proposition was \lq\lq known to the experts\rq\rq but we did not find it in the literature.

\begin{prop}\label{maximal subrank}
The orbit of the unit tensor $\GL(A)\times\GL(B)\times\GL(C)\cdot\unit$ is the set of all maximal subrank tensors.
\end{prop}
\begin{proof}
Let $T$ be a maximal subrank tensor. Then there exist $X\in \End(A)$, $Y\in\End(B)$, and $Z\in\End(C)$ such that  $\unit=(X\otimes Y\otimes Z)\cdot T\in\im(X)\otimes\im(Y)\otimes\im(Z)$. Since $\unit$ is concise, we get that $X,Y,Z$ are invertible.
\end{proof}

Since the subrank and border subrank of a tensor are invariant under multiplying an nonzero scalar, we can consider tensors in the projective space $\PP(A\otimes B\otimes C)$. Let $$\bQ_{\max}:=\overline{\left\{[T]\in\PP (A\otimes B\otimes C):\bQ(T)=n\right\}}$$ be the closure of the set of the maximal border subrank tensors in the projective space. Then $$\overline{\GL(A)\times\GL(B)\times\GL(C)\cdot[\unit]}\subset\bQ_{\max}$$
and we want to find and describe $w\in A\otimes B\otimes C$ such that $\unit+w$ is also of maximal border subrank. 

\begin{thm}[Main Theorem]
    $\dim(\bQ_{\max})\geq(2n^3+3n^2-2n-3)/3.$
\end{thm}

\subsection{Acknowledgements}
The author thanks Harm Derksen for helpful conversations and proposition \ref{maximal subrank}, thanks Fulvio Gesmundo for proposition \ref{generic border subrank}, thanks Austin Conner and Mateusz Michalek for examples in section \ref{nullcone}, thanks JM Landsberg for advising me and his NSF grant AF-1814254 that partially supported this work.


\section{Symmetry group}
Let $\Tilde{G}$ be the group $\GL(A)\times\GL(B)\times\GL(C)$, and let  $\Tilde{\fg}:=\fgl(A)\oplus\fgl(B)\oplus\fgl(C)$ be the Lie algebra of $\Tilde{G}$. For any tensor $T\in A\otimes B\otimes C$, the subgroup of $\Tilde{G}$ preserving $T$ is $$\Tilde{G}_T:=\{g\in \GL(A)\times\GL(B)\times\GL(C):g\cdot T=T\}.$$ Then the Lie algebra of $\Tilde{G}_T$ is $$\Tilde{\fg}_T=\{(x,y,z)\in\fgl(A)\oplus\fgl(B)\oplus\fgl(C):(x,y,z).T=0\}.$$

Since the action of $\GL(A)\times\GL(B)\times\GL(C)$ on $A\otimes B\otimes C$ is not faithful, we will consider the symmetry group of tensors, which is defined as follows: 

Let $\Phi:\GL(A)\times\GL(B)\times\GL(C)\longrightarrow\GL(A\otimes B\otimes C)$ be the natural group action of $\GL(A)\times\GL(B)\times\GL(C)$ on $A\otimes B\otimes C$. Then the kernel of $\Phi$ is $$\ker(\Phi)=\{(\lambda\Id_A,\mu\Id_B,\nu\Id_C):\lambda\mu\nu=1\},$$ whose Lie algebra is the kernel of the differential of $\Phi$ $$\ker(d\Phi)=\{(\lambda\Id_A,\mu\Id_B,\nu\Id_C):\lambda+\mu+\nu=0\}.$$ So the kernel is of dimension two. Let $G:=\GL(A)\times\GL(B)\times\GL(C)/\ker(\Phi)$, and let $\fg$ be its Lie algebra. Then there is a faithful action of $G$ on $A\otimes B\otimes C$. 

For any $T\in A\otimes B\otimes C$, {\it the symmetry group of $T$}, denoted $G_T$, is the stabilizer of $T$ in $G$
$$G_T:=\{g\in G:g\cdot T=T\}.$$ 
The Lie algebra of the symmetry group of $T$ is $$\fg_T=\{X\in\fg:X.T=0\}.$$ 

Note that for any tensor $T$, the subgroup $\Tilde{G}_T$ of $\GL(A)\times\GL(B)\times\GL(C)$ preserving T contains the kernel of $\Phi$ and its Lie algebra $\Tilde{\fg}_{T}$ contains the kernel of $d\Phi$. So the Lie algebra of the symmetry group of $T$ is $\fg_{T}=\Tilde{\fg}_{T}/\ker(d\Phi)$ and the symmetry group of $T$ is $G_T=\Tilde{G}_T/\ker(\Phi)$. 

\subsection{Symmetry group of the unit tensor $\unit$}

Note that the Lie algebra of the group preserving the unit tensor $\unit$ is $$\Tilde{\fg}_{\unit}=\{(x,y,z)\in\fgl(A)\times\fgl(B)\times\fgl(C):(x,y,z).\unit=0\}.$$ Given $(x,y,z)\in\Tilde{\fg}_{\unit}$, we have that 
\begin{align*}
    0&=(x,y,z).\unit=(x,y,z).\sum_{i=1}^n a_i\otimes b_i\otimes c_i \\
    &=\sum_{i=1}^n[x(a_i)\otimes b_i\otimes c_i+a_i\otimes y(b_i)\otimes c_i+a_i\otimes b_i\otimes z(c_i)] \\
    &=\sum_{i=1}^n\left[\sum_{s=1}^n x_{si}a_s\otimes b_i\otimes c_i+a_i\otimes \sum_{t=1}^n y_{ti}b_t\otimes c_i+a_i\otimes b_i\otimes \sum_{u=1}^nz_{ui}c_u\right],
\end{align*}
which implies that $x_{si}=y_{ti}=z_{ui}=0$ for $s,t,u\neq i$ and $x_{ii}+y_{ii}+z_{ii}=0$ for all $i=1,\dots,n$.
So the Lie algebra is $$\Tilde{\fg}_{\unit}=\{(x,y,z):x,y,z\text{ are diagonal $n\times n$ matrices and }x+y+z=0\}.$$
Then the dimension of the subgroup $\Tilde{G}_{\unit}$ preserving the unit tensor is $2n$. Thus the dimension of the symmetry group of the unit tensor is $$\dim G_{\unit}=\dim(\Tilde{G}_{\unit})-2=2n-2.$$
By the proposition \ref{maximal subrank}, the dimension of the set of all maximal subrank tensors is $3n^2-2n$.

The symmetric group on $n$ elements $\fS_n$ can be viewed as a subgroup of $\GL(A)\times\GL(B)\times\GL(C)$ by identifying $\sigma\in\fS_n$ with $(f_\sigma,g_\sigma,h_\sigma)\in\GL(A)\times\GL(B)\times\GL(C)$ where $f_\sigma(a_i)=a_{\sigma(i)}$, $g_\sigma(b_i)=b_{\sigma(i)}$, and $h_\sigma(c_i)=c_{\sigma(i)}$. It is clear that $\fS_n$ and the torus
$$\Tilde{\mathbf{T}}=\left\{(\lambda,\mu,\nu):\lambda,\mu,\nu\text{ are $n\times n$ diagonal matrices and }\lambda\mu\nu={\rm Id}\right\}$$
are both subgroup of $\Tilde{G}_{\unit}$. We claim that $\Tilde{G}_{\unit}=\fS_n\ltimes\Tilde{\mathbf{T}}$. Write the unit tensor $\unit=\sum_{i=1}^n a_i\otimes b_i\otimes c_i$ as a matrix whose entries are in $C$:
$$\left(\begin{array}{ccc}
    c_1 &  & 0 \\
     & \ddots &  \\
    0 &  & c_n
 \end{array}\right)$$
The determinant of this matrix is $c_1c_2\cdots c_n\in S^nC$. For a linear map $T:B^*\longrightarrow A$, we have $T^{\wedge n}:\Lambda^n B^*\longrightarrow\Lambda^nA$ and an element $(f,g)\in\GL(A)\times\GL(B)$ acts on $\Lambda^nT$ by $$(f,g)\cdot T^{\wedge n}=\det(f)\det(g)T^{\wedge n}.$$ The determinant of the matrix presenting $(f,g,h)\cdot \unit$ for $(f,g,h)\in\GL(A)\times\GL(B)\times\GL(C)$ is $$h\cdot(\det(f)\det(g)c_1c_2\cdots c_n)=\det(f)\det(g)h(c_1)h(c_2)\cdots h(c_n).$$

For $(f,g,h)\in\Tilde{G}_{\unit}$, as $\unit=(f,g,h)\cdot\unit$, we have $$c_1c_2\cdots c_n=\det(f)\det(g)h(c_1)h(c_2)\cdots h(c_n)\in S^nC.$$ By the unique factorization of polynomials, there exist $\sigma\in\fS_n$ and $\nu_1,\dots,\nu_n\in\C$ such that $h(c_i)=\nu_i c_{\sigma(i)}$. Similarly, we have $\sigma',\sigma''\in\fS_n$ and $\lambda_1,\dots,\lambda_n,\mu_1,\dots,\mu_n\in\C$ such that $f(a_i)=\lambda_ia_{\sigma'(i)}$ and $g(b_i)=\mu_ib_{\sigma''(i)}$. Since the unit tensor is fixed by $(f,g,h)$, we see that $\sigma=\sigma'=\sigma''$ and $\lambda_i\mu_i\nu_i=1$. So, $$(f,g,h)=\sigma\circ \left(\left(\begin{array}{ccc}
    \lambda_1 &  & 0 \\
     & \ddots &  \\
    0 &  & \lambda_n
 \end{array}\right),\left(\begin{array}{ccc}
    \mu_1 &  & 0 \\
     & \ddots &  \\
    0 &  & \mu_n
 \end{array}\right),\left(\begin{array}{ccc}
    \nu_1 &  & 0 \\
     & \ddots &  \\
    0 &  & \nu_n
 \end{array}\right)\right).$$
 
Given $\sigma\in\fS_n$ and $(\lambda,\mu,\nu)\in\Tilde{\mathbf{T}}$, we have $\sigma(\lambda,\mu,\nu)\sigma^{-1}=(\lambda',\mu',\nu')\in \Tilde{\mathbf{T}}$ where $\lambda'_{ii}=\lambda_{\sigma^{-1}(i)\sigma^{-1}(i)}$, $\mu'_{ii}=\mu_{\sigma^{-1}(i)\sigma^{-1}(i)}$, and $\nu'_{ii}=\nu_{\sigma^{-1}(i)\sigma^{-1}(i)}$.
This gives the semidirect product.

Note that $\Tilde{\mathbf{T}}$ contains the kernel of $\Phi$ and $\fS_n$ has trivial intersection with the kernel of $\Phi$. We can consider the torus $\mathbf{T}:=\Tilde{\mathbf{T}}/\ker(\Phi)$ and see that the symmetry group of the unit tensor is $G_{\unit}=\fS_n\ltimes\mathbf{T}$.



\section{The nullcone $\cN_{G_{\unit}}$} \label{W}
Define the nullcone $\cN_{G_{\unit}}:=\{v\in A\otimes B\otimes C:0\in\overline{G_{\unit}\cdot v}\}$, and let
$$\cC:=\Cone(\unit,\cN_{G_{\unit}})=\left\{\left[\unit+w\right]\in\PP(A\otimes B\otimes C):w\in\cN_{G_{\unit}}\right\}$$ be the projective cone over $\PP\cN_{G_{\unit}}$ with vertex $\left[\unit\right]$.
Then we have $\overline{G\cdot\cC}\subset\bQ_{\max}$.

Note that the nullcone defined by the maximal torus $\mathbf{T}$, $$\cN_{\mathbf{T}}:=\{v\in A\otimes B\otimes C:0\in\overline{\mathbf{T}\cdot v}\}$$ is a subset of the nullcone $\cN_{G_{\unit}}$. By the Hilbert-Mumford criterion, we have $\cN_{G_{\unit}}=G_{\unit}\cdot\cN_{\mathbf{T}}$. Thus we have that $\cN_{G_{\unit}}=\cN_{\mathbf{T}}$ as $\mathbf{T}$ is normal in $G_{\unit}$

Let $x_{ijk}=\alpha_i\otimes\beta_j\otimes\gamma_k\in Sym(A\otimes B\otimes C)^*$. Then the coordinate ring of $A\otimes B\otimes C$ is the polynomial ring $\C[x_{ijk}]$.
The nullcone $\cN_{G_{\unit}}=\cN_{\mathbf{T}}$ is the zero set of all homogeneous polynomials of positive degree that are invariant under the action of $\mathbf{T}$, see e.g., \cite[Lemma 2.4.2]{MR1918599}. Note that $\mathbf{T}$ is a torus, and the monomials in $\C[x_{ijk}]$ are weight vectors of the action of $\mathbf{T}$:
$$\overline{(\lambda,\mu,\nu)}\cdot x_{ijk}=\frac{1}{\lambda_i\mu_j\nu_k}x_{ijk},$$
where $\overline{(\lambda,\mu,\nu)}$ is the image of $(\lambda,\mu,\nu)$ in $\mathbf{T}=\Tilde{\mathbf{T}}/\ker(\Phi)$ and $\lambda_i,\mu_j,\nu_k$ are the $i$-th, $j$-th, and $k$-th diagonal term of $\lambda,\mu,\nu$, respectively.
So the nullcone $\cN_{G_{\unit}}$ is the zero set of monomials that are invariant under the action of $\mathbf{T}$, which implies that the nullcone is the union of linear spaces. 
Let $I$ be the ideal generated by the monomials $$x_{iii},x_{iij}x_{jji},x_{iji}x_{jij},x_{ijj}x_{jii},x_{ijk}x_{jki}x_{kij}$$ for distinct $i,j,k\in\{1,\dots,n\}$. Since those monomials are invariant under the action of $\mathbf{T}$, the nullcone $\cN_{G_{\unit}}$ is contained in the zero set of the ideal $I$. Since no variables in the generators of $I$ are repetitive and the number of the generators is $n+3\binom{n}{2}+2\binom{n}{3}$, the zero set of $I$ is a union of linear spaces with dimension $n^3-(n+3\binom{n}{2}+2\binom{n}{3})$.

Let $W$ be a vector subspace of $A\otimes B\otimes C$ spanned by elements $a_i\otimes b_j\otimes c_k$ satisfying $1\leq i,j,k\leq n$ and at least one of $j,k$ less than $i$. For any $T=\sum_{i,j,k}T_{ijk}a_i\otimes b_j\otimes c_k\in W$, note that $T_{ijk}$ nonzero only when one of $j,k$ is less than $i$, so we have that
\begin{align*}
    &x_{iii}(T)=T_{iii}=0\text{ for all }i=1,\dots,n \\
    &x_{iij}x_{jji}(T)=T_{iij}T_{jji}=0\cdot T_{jji}=0\text{ for all }1\leq i<j\leq n \\
    &x_{iji}x_{jij}(T)=T_{iji}T_{jij}=0\cdot T_{jij}=0\text{ for all }1\leq i<j\leq n \\
    &x_{ijj}x_{jii}(T)=T_{ijj}T_{jii}=0\cdot T_{jii}=0\text{ for all }1\leq i<j\leq n \\
    &x_{ijk}x_{jki}x_{kij}(T)=0\cdot T_{jki}T_{kij}=0\text{ for all }1\leq i<j<k\leq n\text{ or }1\leq i<k<j\leq n.
\end{align*}
Then we have that $W$ is contained in the zero set of $I$. Note that $W$ has dimension $$\dim(W)=n^3-(n^2+(n+2)^2+\cdots+2^2+1^2)=n^3-\frac{n(n+1)(2n+1)}{6},$$ which is equal to the dimension of the zero set of $I$.  Consider $$c(t)=\left(\left(\begin{array}{ccc}
    t^{\lambda_1} &  & 0 \\
     & \ddots &  \\
    0 &  & t^{\lambda_n}
 \end{array}\right),\left(\begin{array}{ccc}
    t^{\mu_1} &  & 0 \\
     & \ddots &  \\
    0 &  & t^{\mu_n}
 \end{array}\right),\left(\begin{array}{ccc}
    t^{\nu_1} &  & 0 \\
     & \ddots &  \\
    0 &  & t^{\nu_n}
 \end{array}\right)\right)\in G_{\unit}$$
where $\lambda_k=2^n-2^{n-k+1}$ and $\mu_k=\nu_k=2^{n-k}-2^{n-1}$ for $k=1,\dots,n$.
For $a_i\otimes b_j\otimes c_k\in W$, we have that 
\begin{align*}
    \lim_{t\to 0}c(t)\cdot(a_i\otimes b_j\otimes c_k)&=\lim_{t\to 0}t^{\lambda_i+\mu_j+\nu_k}(a_i\otimes b_j\otimes c_k) \\
    &=\lim_{t\to 0}t^{2^{n-j}+2^{n-k}-2^{n-i+1}}(a_i\otimes b_j\otimes c_k)=0
\end{align*}
since $2^{n-j}+2^{n-k}\geq2^{n-i+1}+1$. By the definition of the null cone, we have that $W\subset\cN_{G_{\unit}}$. Since the nullcone $\cN_{G_{\unit}}$ is in between two varieties having the same dimension, we can focus on the subset $W$, which is of dimension $$\dim(W)=n^3-\frac{n(n+1)(2n+1)}{6}=\frac{4n^3-3n^2-n}{6}.$$ 
Let $\cC':=\Cone(\unit,W)$ be the projective cone over $\PP W$ with vertex $\left[\unit\right]$. Then its dimension is $$\dim(\cC')=\dim\Cone(\unit,W)=\frac{4n^3-3n^2-n}{6}.$$


\section{The lower bound of the dimension of $\bQ_{\max}$}

Note that $\cC'\subset\cC\subset\bQ_{\max}$ and $\bQ_{\max}$ is invariant under the action of $G$. The dimension of $\overline{G\cdot\cC'}$ give us a lower bound of the dimension of $\bQ_{\max}$. 

\begin{proof}[Proof of the main theorem]
Consider the group preserving $\cC'$: 
\begin{align*}
    G_{\cC'}&:=\{g\in G:g\cdot \cC'\subset\cC'\} \\
    &=\{g\in G:g\cdot[\unit+w]\in\Cone(\unit,W)~\forall w\in W\} \\
    &=\{g\in G:[g\cdot(\unit+w)]\in\Cone(\unit,W)~\forall w\in W\} \\
    &=\{g\in G:g\cdot(\unit+w)=\lambda\unit+w'\text{ for some }\lambda\in\C\setminus\{0\}\text{ and } w'\in W~\forall w\in W\}
\end{align*}
Let $\fg_{\cC'}$ be its Lie algebra. Then  $$\fg_{\cC'}\subset\fg':=\{x\in\fg:x.(\unit+w)\in\langle\unit,W\rangle~\forall w\in W\}.$$ Consider a Lie subalgebra of $\Tilde{\fg}=\fgl(A)\oplus\fgl(B)\oplus\fgl(C)$: $$\Tilde{\fg}':=\{(x,y,z)\in\Tilde{\fg}:(x,y,z).(\unit+w)\in\langle\unit,W\rangle~\forall w\in W\}.$$
Then we have $\fg_{\cC'}\subset\fg'=\Tilde{\fg}'/\ker(d\Phi)$. For $(x,y,z)\in\Tilde{\fg}'$, we have the following equations for all $w=\sum_{i,j,k}w_{ijk}a_i\otimes b_j\otimes c_k\in W$,
\begin{align*}
    (\sigma,\sigma,\sigma)&~~~~ x_{\sigma\sigma}+y_{\sigma\sigma}+z_{\sigma\sigma}+\sum_{i>\sigma}w_{i\sigma\sigma}x_{\sigma i}+\sum_{j<\sigma}w_{\sigma j\sigma}y_{\sigma j}+\sum_{k<\sigma}w_{\sigma\sigma k}z_{\sigma k} \\
    &=x_{11}+y_{11}+z_{11}+\sum_{i>1}w_{i11}x_{1i}\text{ for all }1\leq\sigma\leq n \\
    (\sigma,\rho,\rho)&~~~~ x_{\sigma\rho}+\sum_{i>\rho}w_{i\rho\rho}x_{\sigma i}+\sum_j w_{\sigma j\rho}y_{\rho j}+\sum_k w_{\sigma\rho k}z_{\rho k}=0\text{ for all }\sigma<\rho \\
    (\sigma,\rho,\sigma)&~~~~ y_{\rho\sigma}+\sum_i w_{i\rho\sigma}x_{\sigma i}+\sum_{j<\sigma} w_{\sigma j\sigma}y_{\rho j}+\sum_k w_{\sigma\rho k}z_{\sigma k}=0\text{ for all }\sigma<\rho \\
    (\sigma,\sigma,\rho)&~~~~ z_{\rho\sigma}+\sum_i w_{i\sigma\rho}x_{\sigma i}+\sum_j w_{\sigma j\rho}y_{\sigma j}+\sum_{k<\sigma}w_{\sigma\sigma k}z_{\rho k}=0\text{ for all }\sigma<\rho \\
    (\sigma,\rho,\tau)&~~~~ \sum_i w_{i\rho\tau}x_{\sigma i}+\sum_j w_{\sigma j\tau}y_{\rho j}+\sum_k w_{\sigma \rho k} z_{\tau k}=0\text{ for all }\sigma<\rho,\tau\
\end{align*}
For $\sigma<\rho$, we use induction on $\sigma$ to prove that $x_{\sigma\rho}=y_{\rho\sigma}=z_{\rho\sigma}=0$. For $\sigma=1$ and $\rho>1$, equation $(\sigma,\rho,\rho)$ implies that $x_{1\rho}=0$ since $w_{1jk}=0$ for any $j,k$. Then equations $(\sigma,\rho,\sigma)$ and $(\sigma,\sigma,\rho)$ imply that $y_{\rho 1}=z_{\rho 1}=0$ for all $\rho>1$. Now fix some $\sigma>1$. Assume that $x_{\sigma'\rho}=y_{\rho\sigma'}=z_{\rho\sigma'}=0$ for all $\sigma'<\sigma$ and $\rho>\sigma'$. For all $\rho>\sigma$, by the induction hypothesis, we have the equations 
\begin{align*}
    (\sigma,\rho,\rho)&~~~~ x_{\sigma\rho}=x_{\sigma\rho}+\sum_{j<\sigma} w_{\sigma j\rho}y_{\rho j}+\sum_{k<\sigma} w_{\sigma\rho k}z_{\rho k}=0 \\
    (\sigma,\rho,\sigma)&~~~~ y_{\rho\sigma}+\sum_{\rho>i>\sigma} w_{i\rho\sigma}x_{\sigma i}=y_{\rho\sigma}+\sum_{\rho>i>\sigma} w_{i\rho\sigma}x_{\sigma i}+\sum_{k<\sigma} w_{\sigma\rho k}z_{\sigma k}=0\\
    (\sigma,\sigma,\rho)&~~~~ z_{\rho\sigma}+\sum_{\rho>i>\sigma} w_{i\sigma\rho}x_{\sigma i}=z_{\rho\sigma}+\sum_{\rho>i>\sigma} w_{i\sigma\rho}x_{\sigma i}+\sum_{j<\sigma} w_{\sigma j\rho}y_{\sigma j}=0.
\end{align*}
Then we have that
$$(x,y,z)=\left(\left(\begin{array}{cccc}
    x_{11} & 0 & \cdots & 0 \\
    \vdots & \ddots & \ddots & \vdots \\
    \vdots &  & \ddots & 0 \\
    x_{n1} & \cdots & \cdots & x_{nn}
 \end{array}\right),\left(\begin{array}{cccc}
    y_{11} & \cdots & \cdots & y_{1n} \\
    0 & \ddots & & \vdots \\
    \vdots & \ddots & \ddots & \vdots \\
    0 & \cdots & 0 & y_{nn}
 \end{array}\right),\left(\begin{array}{cccc}
    z_{11} & \cdots & \cdots & z_{1n} \\
    0 & \ddots & & \vdots \\
    \vdots & \ddots & \ddots & \vdots \\
    0 & \cdots & 0 & z_{nn}
 \end{array}\right)\right)$$
satisfying $$(\sigma,\sigma,\sigma)~~~~ x_{\sigma\sigma}+y_{\sigma\sigma}+z_{\sigma\sigma}=x_{11}+y_{11}+z_{11}\text{ for all }1<\sigma\leq n.$$ 
Therefore we have $\fg_{\cC'}=\fg'$ and $\fg_{\cC'}$ is of dimension $3(n+1)n/2-(n-1)-2=(3n^2+n-2)/2$. So $\dim(G_{\cC'})=\dim(\fg_{\cC'})=(3n^2+n-2)/2$. Then the dimension of $\overline{G\cdot\cC'}$ is
\begin{align*}
    \dim(G)&-\dim(G_{\cC'})+\dim(\cC') \\
    &=(3n^2-2)-(3n^2+n-2)/2+(4n^3-3n^2-n)/6 \\
    &=(2n^3+3n^2-2n-3)/3,
\end{align*}
which is a lower bound of the dimension of $\bQ_{\max}$.

\end{proof}


\section{Generic border subrank is not maximal}
\begin{prop} \label{generic border subrank}
    The generic border subrank is at most $n-1$.
\end{prop}
\begin{proof}
    Let $v\in A\otimes B\otimes C$, and let $S=\SL(A)\times\SL(B)\times\SL(C)$. Since the unit tensor $\unit$ is semistable under the action of $S$, there is a homogeneous polynomial $P_1$ that is invariant under the action of $S$ and $P_1(\unit)\neq 0$. Let $P_2$ be another homogeneous polynomial that is invariant under the action of $S$ and is algebraically independent with $P_1$. Without loss of generality, we may assume that the degree of $P_1$ and $P_2$ are the same, say $d$. Let $\alpha_i=P_i(v)$ and $\beta_i=P_i(\unit)$ for $i=1,2$. Then $\alpha_2P_1-\alpha_1P_2$ is also invariant under the action of $S$. For $(f,g,h)\in\GL(A)\times\GL(B)\times\GL(C)$, there are nonzero $\lambda,\mu,\nu\in\C$ such that $(f/\lambda,g/\mu,h/\nu)\in S$, and so
    \begin{align*}
        (\alpha_2P_1-\alpha_1P_2)((f,g,h)\cdot v)&=(\alpha_2P_1-\alpha_1P_2)(\lambda\mu\nu(f/\lambda,g/\mu,h/\nu)\cdot v) \\
        &=(\lambda\mu\nu)^d(\alpha_2P_1-\alpha_1P_2)((f/\lambda,g/\mu,h/\nu)\cdot v)\\
        &=(\lambda\mu\nu)^d(\alpha_2P_1-\alpha_1P_2)(v) \\
        &=(\lambda\mu\nu)^d(\alpha_2\alpha_1-\alpha_1\alpha_2)=0.
    \end{align*}
    Thus, the orbit closure of $v$, $\overline{G\cdot v}$, is contained in the zero set of $\alpha_2P_1-\alpha_1P_2$.
    Note that $$(\alpha_2P_1-\alpha_1P_2)(\unit)=\alpha_2\beta_1-\alpha_1\beta_2=(\beta_1 P_2-\beta_2 P_1)(v).$$  If $v$ is not in the zero set of $\beta_1 P_2-\beta_2 P_1$, then $(\alpha_2P_1-\alpha_1P_2)(\unit)$ is nonzero, which implies that the unit tensor $\unit$ is not in the orbit closure $\overline{G\cdot v}$ and so the subrank of $v$ is not maximal. We proved that the generic border subrank is not maximal.
\end{proof}


\section{Find the nullcone $\cN_{G_{\unit}}$} \label{nullcone}
Recall that $$W=\langle a_i\otimes b_j\otimes c_k:\text{at least one of }j,k\text{ less than }i\rangle$$ is of dimension $(4n^3-3n^2-n)/6$. Similarly, we can define $$W':=\langle a_i\otimes b_j\otimes c_k:\text{at least one of }i,k\text{ less than }j\rangle$$ and $$W'':=\langle a_i\otimes b_j\otimes c_k:\text{at least one of }i,j\text{ less than }k\rangle.$$ Remark that $W'$ and $W''$ are both subsets of the nullcone $\cN_{G_{\unit}}$. For $\sigma\in\fS_n$, define $$\sigma(W):=\langle a_{\sigma(i)}\otimes b_{\sigma(j)}\otimes c_{\sigma(k)}:\text{at least one of }j,k\text{ less than }i\rangle.$$ $\sigma(W')$ and $\sigma(W'')$ are defined similarly. Then we have $$\bigcup_{\sigma\in\fS_n}(\sigma(W)\cup\sigma(W')\cup\sigma(W''))\subset\cN_{G_{\unit}}.$$
However this is not the whole nullcone and in fact that the nullcone is not equidimensional.

\begin{ex}
    Consider a one-parameter subgroup $$c(t)=\left(\left(\begin{array}{ccc}
    t^5 &  & 0 \\
     & 1 &  \\
    0 &  & t^2
 \end{array}\right),\left(\begin{array}{ccc}
    1 &  & 0 \\
     & t &  \\
    0 &  & t^{-3}
 \end{array}\right),\left(\begin{array}{ccc}
    t^{-5} &  & 0 \\
     & t^{-1} &  \\
    0 &  & t
 \end{array}\right)\right)\in G_{M^{\oplus 3}_{\langle1\rangle}}.$$ Let $U$ be a subspace of $A\otimes B\otimes C$ generated by $a_i\otimes b_j\otimes c_k$ where $\lim_{t\to 0}c(t)\cdot a_i\otimes b_j\otimes c_k=0$. Then $U$ is a subset of the nullcone $\cN_{G_{M^{\oplus 3}_{\langle1\rangle}}}$ and is of dimension $13=\dim W$, but $U$ is not any of $\sigma(W),\sigma(W'),\sigma(W'')$.
\end{ex}

\begin{ex}
Consider a one-parameter subgroup $$c(t)=\left(\left(\begin{array}{ccc}
    t^{-2} &  & 0 \\
     & t^{-1} &  \\
    0 &  & 1
 \end{array}\right),\left(\begin{array}{ccc}
    t^3 &  & 0 \\
     & t^{-2} &  \\
    0 &  & 1
 \end{array}\right),\left(\begin{array}{ccc}
    t^{-1} &  & 0 \\
     & t^3 &  \\
    0 &  & 1
 \end{array}\right)\right)\in G_{M^{\oplus 3}_{\langle1\rangle}}.$$ Let $U$ be a subspace of $A\otimes B\otimes C$ generated by $a_i\otimes b_j\otimes c_k$ where $\lim_{t\to 0}c(t)\cdot a_i\otimes b_j\otimes c_k=0$. Then $U$ is a subset of the nullcone $\cN_{G_{M^{\oplus 3}_{\langle1\rangle}}}$ and is of dimension $12<13=\dim W$. By considering the monomials listed in section 3 and these three monomials $x_{123}x_{211}x_{332},$ $x_{231}x_{322}x_{113},$ and $x_{132}x_{321}x_{213}$, we can see that $U$ is a maximal component in the nullcone $\cN_{G_{M^{\oplus 3}_{\langle1\rangle}}}$. So the null cone $\cN_{G_{M^{\oplus 3}_{\langle1\rangle}}}$ is not equidimensional.
\end{ex}


\section{Tight and maximal border subrank tensors}
\begin{defi}
    The {\it support} of a tensor $T\in A\otimes B\otimes C$, denote $supp(T)$, is $$supp(T):=\{(i,j,k):T_{ijk}\neq 0\}$$ where $T=\sum_{ijk}T_{ijk}a_i\otimes b_j\otimes c_k$. A tensor $T$ is called {\it tight} if there are injective function $\tau_A,\tau_B,\tau_C:\{1,\dots,n\}\rightarrow\Z$ such that $\tau_A(i)+\tau_B(j)+\tau_C(k)=0$ for all $(i,j,k)\in supp(T)$.
\end{defi}

To use the Laser method, we need ``blocked tight'' tensors, and tight tensors are blocked tight. Tensors that have been used for the Laser method are actually tight, not just blocked tight. We show that there are tight and maximal border subrank tensors.

Consider the linear subspace $$U:=\langle a_i\otimes b_j\otimes c_k:2i=j+k\text{ and at least one of }j,k\text{ less than }i\rangle\subset W\subset\cN_{G_{\unit}}.$$ Then $\Cone(\unit,U)\subset\Cone(\unit,W)\subset\Cone(\unit,\cN_{G_{\unit}})\subset\bQ_{\max}$. For any tensor $T\in A\otimes B\otimes C$ with $[T]\in\Cone(\unit,U)$, its support is
\begin{align*}
    supp(T)&=\{(i,i,i):1\leq i\leq n\}\cup\{(i,j,k):2i=j+k\text{ and at least one of }j,k\text{ less than }i\} \\
    &=\{(i,i,i):1\leq i\leq n\}\cup\{(i,j,k):2i=j+k\text{ and }j,k\neq i\} \\
    &=\{(i,j,k):2i=j+k\}
\end{align*}
Define the injective functions $\tau_A,\tau_B,\tau_C:\{1,\dots,n\}\rightarrow\Z$ by $\tau_A(i)=3-2i$, $\tau_B(j)=j$, and $\tau_C(k)=k-3$. Since $\tau_A(i)+\tau_B(j)+\tau_C(k)=0$ for all $(i,j,k)\in supp(T)$, this shows that $T$ is a tight tensor. Therefore, $\Cone(\unit,U)$ contains tight and maximal border subrank tensors.



\bibliographystyle{alpha}
\bibliography{reference}

\end{document}